\documentclass[10pt,twoside,english,1p]{elsarticle}

\usepackage[T1]{fontenc}
\usepackage[latin9]{inputenc}
\usepackage{geometry}
\usepackage{color}
\usepackage{float}
\usepackage{units}
\usepackage{amsmath}
\usepackage{amsthm}
\usepackage{graphicx}
\usepackage{setspace}
\usepackage{esint}
\makeatletter

\usepackage{pgf,tikz,pgfplots,adjustbox, mathrsfs, lineno, booktabs}

\pgfplotsset{compat=1.15}
\usetikzlibrary{arrows}
\usetikzlibrary{external}
\usepackage{framed}
\usepackage{multirow}

\usepackage{latexsym}

\usepackage{url}
\usepackage{xcolor, soul}
\definecolor{newcolor}{rgb}{.8,.349,.1}
\usepackage{xcolor-material}
\colorlet{R2}{MaterialTeal900}
\colorlet{R1}{MaterialIndigo900}

\colorlet{B2}{MaterialTeal100}
\colorlet{B1}{MaterialIndigo100}

\journal{Computers \& Structures}

\definecolor{teal}{RGB}{0,77,64}
\definecolor{blue}{RGB}{26,35,126}
\usepackage{nicefrac}

\expandafter\def\expandafter\normalsize\expandafter{%
  \normalsize  
  \setlength\abovedisplayskip{0.0001ex}
  \setlength\belowdisplayskip{1ex}
  \setlength\abovedisplayshortskip{1ex}
  \setlength\belowdisplayshortskip{1ex}
}

\@ifundefined{showcaptionsetup}{}{%
 \PassOptionsToPackage{caption=false}{subfig}}
\usepackage{subfig}
\makeatother

\usepackage{babel}
\begin{document}

\begin{frontmatter}{}

\title{\textbf{A Finite Volume scheme for the solution of discontinuous
magnetic field distributions on non-orthogonal meshes}}

\author[b]{Augusto Riedinger}

\author[b]{Martín Saravia\corref{cor1}}

\author[b]{José Ramírez}

\ead{msaravia@conicet.gov.ar}

\cortext[cor1]{Corresponding author}

\address[b]{Centro de Investigación en Mecánica Teórica y Aplicada, Universidad
Tecnológica Nacional, Facultad Regional Bahía Blanca, CONICET. 11
de Abril 461, Bahía Blanca, Argentina.}
\begin{abstract}
We present a Finite Volume formulation for determining discontinuous
distributions of magnetic fields within non-orthogonal and non-uniform
meshes. The numerical approach is based on the discretization of the
vector potential variant of the equations governing static magnetic
field distribution in magnetized, permeable and current carrying media.
After outlining the derivation of the magnetostatic balance equations
and its associated boundary conditions, we propose a cell\textendash centered
Finite Volume framework for spatial discretization and a Block Gauss\textendash Seidel
multi-region scheme for solution. We discuss the structure of the
solver, emphasizing its effectiveness and addressing stabilization
and correction techniques to enhance computational robustness. We
validate the accuracy and efficacy of the approach through numerical
experiments and comparisons with the Finite Element method.
\end{abstract}

\end{frontmatter}{}


\newpage{}

\section{Introduction}

Magnetic fields have captivated researchers since their discovery,
and in recent years, there has been a surge in their engineering and
industrial applications. Consequently, the field of magnetostatics
has experienced significant advancements. However, addressing the
intricate geometries encountered in real-world applications has always
been a major challenge due to the inherent complexity of Maxwell's
equations\textemdash the fundamental partial differential equations
governing the dynamics of magnetic fields \citep{maxwell1865viii}.
Thus, the application of numerical methods for predicting magnetic
fields has become indispensable in tackling this intricate problem
and seeking practical solutions.x

One of the primary consequences of the technological advances in fabrication
techniques is the emergence of designs with complex geometries. Although
they have proven to enhance the efficiency and appearance of modern
devices, the computational simulation of the physical phenomena in
which they play a role has become increasingly challenging.

In the landscape of Computational Electromagnetics (CEM), researchers
have made remarkable strides since its inception \citep{steele2012numerical,sykulski2012computational,lowther2012computer}..
Pioneering work by T. K. Sarkar et al. laid the groundwork for numerical
schemes aimed at solving linear matrix equations, particularly in
the context of electromagnetic scattering and radiation problems.
Building upon these foundations, subsequent years have witnessed a
surge in advancements, as evidenced by the seminal contributions of
Salon et al \citep{salon1999numerical}. Further enriching the field,
applied investigations have delved into diverse areas, with notable
reviews such as Niu et al.'s comprehensive exploration of electromagnetic
modeling providing valuable insights \citep{niu2020numerical}. Despite
these advancements, the literature addressing the numerical modeling
of electromagnetic fields on complex geometries remains relatively
sparse.

The Finite Element Method (FEM) has traditionally been employed to
numerically solve electromagnetic problems \citep{cardoso2016electromagnetics}.
Early formulations were pioneered by J. H. Coggon et al. \citep{coggon1971electromagnetic}
and Sadiku \citep{sadiku1989simple,sadiku2018numerical}, who provided
a foundational introduction to the FEM applied to electromagnetics.
Subsequent advances were made by J. Jin et al. \citep{jin2015finite}
and A. C. Polycarpou et al. \citep{polycarpou2022introduction}. A
recent review addresing the use of the FEM to model electromagnetic
problems was presented by M. Augustyniak et al. \citep{augustyniak2016finite}.

In contrast, the application of the Finite Volume Method (FVM) to
electromagnetics took more time to develop. It began with C.-D. Munz
et al. \citep{munz2000finite}, who developed a time domain approach
to solve Maxwell equations, later extended by T. S. Chung et al. \citep{chung2001finite}
to general 3D polyhedral domains with discontinuous physical coefficients.
Later, Y. Liu et al. \citep{liu2006spectral} investigated spectral
FVM on unstructured grids, presenting a method computationally more
effective than the conventional structured methods. In parallel, several
applied studies emerged, such as the solution of automotive electromagnetic
compatibility problems using a hybrid finite difference/Finite Volume
method by X. Ferrieres et al. \citep{ferrieres2004application}, and
techniques developed by E. Haber et al. \citep{haber2014multiscale}
and L. A. C. Mata et al. \citep{caudillo2017oversampling} to reduce
the computational cost of solving the Maxwell equations at low frequencies,
with direct applications to geophysical settings. However, the aforementioned
papers encounter limitations in handling highly non-orthogonal grids
and discontinuous magnetic properties across interfaces in domains
composed of media of different physical characteristics.

In 2021, Saravia \citep{saravia2021finite} introduced a FVM framework
to solve discontinuous distributions of magnetic fields, employing
a multi-region approach implemented within the OpenFOAM library. The
author presented a formulation written in terms of the magnetic vector
potential to determine the distribution of magnetic fields in interacting
permeable, permanently magnetized and current carrying media. They
benchmarked the accuracy of the FVM multi-region magnetostatic approach
against the FEM. Subsequently, Riedinger et al. \citep{riedinger2023single}
expanded the framework for single-region schemes, demonstrating that
the magnetic field could be accurately captured with a continuous
method using boundary layer meshing. Although the mentioned methods
have proven to be effective to solve magnetostatic problems, they
have not been tested in non-orthogonal meshes.

The literature concerning the numerical treatment of discontinuous
distributions of magnetic fields using the Finite Volume Method is
sparse. In particular, the application of this method to obtain solutions
to Maxwell's equations on non orthogonal meshes is even scarser. Nevertheless,
the numerous instances of FVM application in other areas of physics
which involve complex geometries and non-orthogonal meshes, such as
fluid mechanics and heat transfer, suggest that the FVM may still
be an appealing option.

In this context, we present a Finite Volume formulation tailored to
analyze discontinuous distributions of magnetic fields within non-orthogonal
and non-uniform meshes. Our numerical approach centers on discretizing
the vector potential variant of the equation governing static magnetic
field distribution in magnetized, permeable, and current-carrying
media. We provide an overview of the derivation of the magnetostatic
balance equations and their associated boundary conditions. Subsequently,
we introduce a cell-centered Finite Volume framework for spatial discretization
and propose a Block Gauss\textendash Seidel multi-region scheme for
solution. Emphasizing the structure of the solver, we highlight its
effectiveness and discuss stabilization and limiting techniques aimed
at enhancing computational robustness. Through numerical experiments
and comparisons with the Finite Element Method, we validate the accuracy
and efficacy of the formulation.

\section{Magnetostatic formulation\label{sec:Fundamental-magnetostatic-equati}}

\subsection{Conservation laws}

We begin assuming that currents are not changing with time. Then,
Maxwell's equations may be separated into two sets: one for the electric
field and another for the magnetic field. Under this conditions, magnetic
fields are static, and the conservation laws reduce to
\begin{equation}
\nabla\cdot\mathbf{B}=0\label{eq:gauss-law}
\end{equation}

and
\begin{equation}
\nabla\times\mathbf{B}=\mu_{0}\,\mathbf{J}_{f}\label{eq:ampere-law}
\end{equation}
for the magnetic field $\mathbf{B}$,and 
\begin{equation}
\nabla\cdot\mathbf{J}_{f}=0\label{eq:free-current-conservation-1}
\end{equation}

for the free current $\mathbf{J}_{f}$ \citep{jackson1999classical,griffiths2021introduction}.
Equation \eqref{eq:gauss-law} is Gauss' law for magnetism and \eqref{eq:ampere-law}
is Ampere's law.

The numerical solution to these equations within the context of permeable,
magnetized, and current-carrying media necessitates formulating equilibrium
equations and associated boundary conditions governing the behavior
of the magnetic field. In the subsequent sections, we outline the
derivation process.

\subsection{Ampere's law}

\subsubsection{Permeable media and current carrying media}

The distribution of the magnetic fields within a permeable domain
is governed by a variant of Ampere's law given by
\begin{equation}
\nabla\times\mathbf{B}_{p}=\mu_{0}\,\left(\mathbf{J}_{f}+\mathbf{J}_{i}\right),\label{eq:nabla-times-bp}
\end{equation}
being $\mu_{0}\approx4\,\pi\times10^{-7}\,\left[\nicefrac{\text{H}}{\text{m}}\right]$
the magnetic permeability of free space, and $\mathbf{J}_{i}$ the
so-called induced current \citep{saravia2021finite}.

Interpreting $\mathbf{J}_{i}$ as the outcome of an induced magnetization
$\mathbf{M}_{i}$ allows the following relation to hold:
\begin{equation}
\mathbf{J}_{i}=\nabla\times\mathbf{M}_{i}.\label{eq:induced-current}
\end{equation}

If we rewrite Eq. \eqref{eq:nabla-times-bp} in the form
\begin{equation}
\nabla\times\mathbf{B}_{p}=\mu_{0}\,\left(\mathbf{J}_{f}+\nabla\times\mathbf{M}_{i}\right),\label{eq:nabla-times-bp-2}
\end{equation}

and re-arrange terms, then 
\begin{equation}
\nabla\times\left(\frac{1}{\mu_{0}}\mathbf{B}_{p}-\mathbf{M}_{i}\right)=\mathbf{J}_{f}.\label{eq:free-current}
\end{equation}

To simplify this expression, it is customary to introduce the auxiliary
field $\mathbf{H}$, defined as
\begin{equation}
\mathbf{H}\equiv\frac{1}{\mu_{0}}\mathbf{B}_{p}-\mathbf{M}_{i}+\nabla\phi,\label{eq:H-definition}
\end{equation}
being $\phi$ is an arbitrary scalar field. It is noteworthy that,
since the curl of a gradient is identically zero, we are free to choose
any $\phi$ and still satisfy Eq. \eqref{eq:free-current}. In particular,
choosing $\text{\ensuremath{\nabla}}\phi=0$ simplifies Equation \eqref{eq:H-definition}
to
\begin{equation}
\mathbf{H}=\frac{1}{\mu_{0}}\mathbf{B}_{p}-\mathbf{M}_{i}.\label{eq:H-simplified}
\end{equation}

Then
\begin{equation}
\nabla\times\mathbf{H}=\mathbf{J}_{f}.\label{eq:free-current-simplified}
\end{equation}

For simplicity, we assume a linear the constitutive law \citep{sommerfeld2013electrodynamics},
so
\begin{equation}
\mathbf{B}_{p}=\mu\mathbf{H},\label{eq:constitutive-law-magnetic-fields}
\end{equation}
being $\mu$ the absolute magnetic permeability of the medium. Expanding
the equation above using the result from Equation \eqref{eq:H-simplified}
gives
\begin{equation}
\mathbf{B}_{p}=\mu\left(\frac{1}{\mu_{0}}\mathbf{B}_{p}-\mathbf{M}_{i}\right).\label{eq:constitutive-law-expanded}
\end{equation}

From this, we can derive an expression for the magnetization $\mathbf{M}_{i}$
as a function of the magnetic field as
\begin{equation}
\mathbf{M}_{i}=\left(\frac{\mu-\mu_{0}}{\mu\,\mu_{0}}\right)\mathbf{B}_{p}=\chi\mathbf{B}_{p},\label{eq:magnetization}
\end{equation}

being $\chi$ the normalized magnetic susceptibility
\begin{equation}
\chi=\frac{\mu-\mu_{0}}{\mu\mu_{0}}=\frac{\mu_{r}-1}{\mu_{r}\mu_{0}},\label{eq:chi}
\end{equation}

and $\mu_{r}=\nicefrac{\mu}{\mu_{0}}$ is the relative permeability.

Substituting Equations \eqref{eq:free-current-simplified} and \eqref{eq:magnetization}
into Equation \eqref{eq:nabla-times-bp-2}, an expanded expression
for the Ampere's law in permeable media carrying a free current can
be obtained:
\begin{equation}
\nabla\times\mathbf{B}_{p}=\mu_{0}\left[\nabla\times\mathbf{H}+\nabla\times\left(\chi\:\mathbf{B}_{p}\right)\right].\label{eq:ampere-permeable-media}
\end{equation}

\subsubsection{Permanently magnetized media}

The distribution of the magnetic field in permanently magnetized media,
from now on denoted as $\mathbf{B}_{m}$, is governed by the following
variant of Ampere's law:
\begin{equation}
\nabla\times\mathbf{B}_{m}=\mu_{0}\:\mathbf{J}_{b}.\label{eq:ampere-law-magnetized-media}
\end{equation}

The so-called bound current $\mathbf{J}_{b}$ is the result of the
permanent magnetization; so, both fields relate through \citep{griffiths2021introduction}
\begin{equation}
\mathbf{J}_{b}=\nabla\times\mathbf{M},\label{eq:bound-current}
\end{equation}

being $\mathbf{M}$ the magnetization vector. We can expand Equation
\eqref{eq:ampere-law-magnetized-media} incorporating Equation \eqref{eq:bound-current}
to obtain
\begin{equation}
\nabla\times\mathbf{B}_{m}=\mu_{0}\:\nabla\times\mathbf{M},\label{eq:ampere-law-magnetized-media-expanded}
\end{equation}

which represents the final form of the Ampere's law for magnetized
media.

\subsubsection{Combined form of the Ampere's law}

In systems consisting of both permeable, magnetized and current-carrying
bodies, the resulting magnetic field $\mathbf{B}$ comprises contributions
from both $\mathbf{B}_{p}$ and $\mathbf{B}_{m}$. Then, Ampere's
law yields
\begin{equation}
\nabla\times\mathbf{B}=\nabla\times\mathbf{B}_{p}+\nabla\times\mathbf{B}_{m}=\mu_{0}\left(\mathbf{J}_{f}+\mathbf{J}_{i}+\mathbf{J}_{b}\right).\label{eq:combined-ampere-law}
\end{equation}

Injecting the expressions for the induced and bound currents we obtain
\begin{equation}
\nabla\times\mathbf{B}=\mu_{0}\left[\mathbf{J}_{f}+\nabla\times\left(\chi\:\mathbf{B}_{p}\right)+\nabla\times\mathbf{M}\right].\label{eq:combined-ampere-law-expanded}
\end{equation}

Now, after assuming the relative permeability of the magnetized medium
satisfies $\mathbf{B}_{p}=\mathbf{B}$, we can devise the generalized
Ampere's law for the entire domain as
\begin{equation}
\nabla\times\mathbf{B}=\mu_{0}\left[\mathbf{J}_{f}+\nabla\times\left(\chi\:\mathbf{B}\right)+\nabla\times\mathbf{M}\right].\label{eq:combined-ampere-law-final}
\end{equation}

This equation is valid only if the permeable media is not permanently
magnetized and if the magnetized material has relative permeability
equal to 1.

\subsection{Vector potential formulation\label{sec:Vector-potential-formulation}}

By virtue of Gauss law, the magnetic field is solenoidal; this means
it can be expressed as the curl of an auxiliary vector field $\mathbf{A}$,
normally called the magnetic vector potential, such that
\begin{equation}
\mathbf{B}=\nabla\times\mathbf{A}.\label{eq:A-definition}
\end{equation}

Choosing Coulomb's gauging, i.e. $\nabla\cdot\mathbf{A}=0$ \citep{jackson1999classical},
allows writing Equation \eqref{eq:combined-ampere-law-final} as
\begin{equation}
\nabla^{2}\mathbf{A}=-\mu_{0}\:\left[\mathbf{J}_{f}+\nabla\times\left(\chi\:\mathbf{B}\right)+\nabla\times\mathbf{M}\right].\label{eq:A-3}
\end{equation}

After expanding the curl in the second term in the RHS, the following
expression for the vector potential balance equation can be obtained:
\citep{saravia2021finite}
\begin{equation}
\nabla^{2}\mathbf{A}=\frac{-\mu_{0}}{1-\chi\:\mu_{0}}\:\left[\mathbf{J}_{f}+\left(\nabla\chi\right)\times\left(\nabla\times\mathbf{A}\right)+\nabla\times\mathbf{M}\right].\label{eq:A-7}
\end{equation}

Although this equation is in its own right a balance law for magnetostatics,
its numerical discretization using the FVM leads to non\textendash conservative
scheme due to the presence of the second term in the RHS. However,
exploiting identities \ref{eq:I12}, \ref{eq:I10} and \ref{eq:I11}
while defining $\nabla\widetilde{\mathbf{A}}=\left(\nabla\mathbf{A}-\nabla\mathbf{A}^{T}\right)$,
the following conservative variant can be developed: \citep{saravia2021finite}
\begin{equation}
\nabla^{2}\mathbf{A}=-\mu_{0}\left[\mathbf{J}_{f}-\nabla\cdot\left(\chi\nabla\widetilde{\mathbf{A}}\right)+\nabla\times\mathbf{M}\right].\label{eq:A-15}
\end{equation}

\subsection{Interface boundary condition}

A crucial aspect of the current formulation is that boundary conditions
must be applied to the vector potential at the interfaces between
the various media comprising the magnetostatic system to ensure that
magnetic field discontinuity is calculated properly. As demonstrated
in \citep{saravia2021finite}, the vector potential and the normal
component of the magnetic field remain continuous across the interface.
So, being $\mathbf{e}_{n}$ the interface unit normal vector, the
jumps in the vector potential and in the magnetic field across the
interface are
\begin{equation}
\Delta\mathbf{A}=\mathbf{0},\label{eq:bc1}
\end{equation}
and 
\begin{equation}
\Delta\mathbf{B}\cdot\mathbf{e}_{n}=0.\label{eq:bc2}
\end{equation}

Conversely, the tangential component of the vector potential experiences
a discontinuity at the interface, and thus the magnitude of its jump
is 
\begin{equation}
\Delta\mathbf{B}=\mu_{0}\,\mathbf{K}\times\mathbf{e}_{n},\label{eq:bc3}
\end{equation}
where $\mathbf{K}$ is a generalized surface current defined as 
\begin{equation}
\mathbf{K}=-\mathbf{K}_{f}+\Delta\left(\chi\mathbf{B}+\mathbf{M}\right)\times\mathbf{e}_{n}.\label{eq:bc4}
\end{equation}
Therefore, the existence of this interface current leads to a discontinuity
in the tangential component of the magnetic field, which in turn leads
to a discontinuity in the normal gradient of the vector potential
of magnitude \citep{saravia2021finite}
\begin{equation}
\frac{\partial\Delta\mathbf{A}}{\partial e_{n}}=-\mu_{0}\mathbf{K}.\label{eq:bc-7}
\end{equation}

The equations above provide the relation between of the magnetic field
and the vector potential on both sides of interfaces that must be
fulfilled in order to satisfy the magnetostatic equilibrium equations.

\section{Numerical scheme\label{sec:Finite-volume-discretization}}

\subsection{Finite Volume formulation}

The primary objective of this paper is to introduce a Finite Volume
scheme for the numerical solution of Equation \eqref{eq:A-15} within
a magnetostatic system comprised of subsystems consisting of permeable,
magnetized, and current-carrying media. The scheme involves dividing
each subsystem into a finite set of arbitrarily irregular polyhedral
cells. These cells are fixed in space, forming a volume over which
each operator of the magnetostatic differential equations and corresponding
interface and boundary conditions are discretized. The results of
this discretization is a linear system of algebraic equations that
are solved in a sequential manner.

We propose to discretize and integrate the balance equations in a
region-wise approach. So, for each region encompassing a subsystem
we define a control volume $V_{i}$ and integrate Equation \eqref{eq:A-15}
as

\begin{equation}
\intop_{V_{i}}\text{\ensuremath{\nabla^{2}}}\mathbf{A}\,\text{d}V=-\mu_{0}\intop_{V_{i}}\left[\mathbf{J}_{f}-\text{\text{\ensuremath{\nabla\cdot}}}\left(\chi\text{\ensuremath{\nabla}}\widetilde{\mathbf{A}}\right)+\text{\ensuremath{\nabla\times}}\mathbf{M}\right]\,\text{d}V.\label{eq:dis1}
\end{equation}

At the end of the discretization process we obtain $n$ linear systems
of equations of the form
\begin{equation}
\mathbf{S}_{i}\mathbf{A}_{i}=\mathbf{F}_{i},
\end{equation}
being $\mathbf{S}_{i}$ the system matrix for region $i$ and $\mathbf{F}_{i}$
the forcing vector resulting from the explicit discretization of the
RHS of Eq. \eqref{eq:dis1} and other explicit terms associated to
the interface boundary conditions and non-orthogonal corrections.

We adopt a segregated solution strategy; so each of the subsystems
is solved independently and sequentially, and the interaction between
them is accounted for by the numerical discretization of the interface
boundary conditions.

\subsubsection{Laplacian discretization}

The Laplacian operator is discretized with the help of Gauss' divergence
theorem as

\begin{equation}
\intop_{V_{i}}\text{\ensuremath{\nabla^{2}}}\mathbf{A}\,\text{d}V=\intop_{V_{i}}\text{\ensuremath{\nabla\cdot}}\left(\text{\ensuremath{\nabla}}\mathbf{A}\right)\,\text{d}V\approx\sum\mathbf{s}_{f}\cdot\text{\ensuremath{\nabla}}\mathbf{A}_{f},
\end{equation}
where $\mathbf{s}_{f}$ is the surface vector, and $\mathbf{A}_{f}$
is the magnetic vector potential at face $f$. The last equation can
also be written in terms of the surface normal gradient as
\[
\sum\mathbf{s}_{f}\cdot\text{\ensuremath{\nabla}}\mathbf{A}_{f}=\sum\left|\mathbf{s}_{f}\right|\hat{\mathbf{n}}\cdot\text{\ensuremath{\nabla}}\mathbf{A}_{f}.
\]

The normal gradient of the magnetic vector potential at the faces
\textendash or surface normal gradient\textendash{} is required to
evaluate the Laplacian. For orthogonal meshes its numerical evaluation
is based on a central difference of cell values on each side of the
face,
\begin{equation}
\hat{\mathbf{n}}\cdot\text{\ensuremath{\nabla}}\mathbf{A}_{f}=\hat{\mathbf{r}}\cdot\text{\ensuremath{\nabla}}\mathbf{A}_{f}=\frac{\mathbf{A}_{C}-\mathbf{A}_{E}}{\left|\mathbf{r}\right|},\label{eq:rGrad}
\end{equation}
being $\mathbf{r}$ the relative position vector between centers of
cells $C$ and $E$, which share face $f$.

\begin{figure}[H]
\begin{centering}
\includegraphics{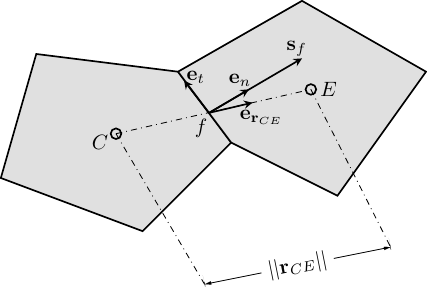}
\par\end{centering}
\caption{Depiction of neighboring elements in a non-orthogonal mesh. \label{fig:non-orthogonal-mesh}}
\end{figure}

This calculation maintains second-order accuracy when the vector connecting
the cell centers is orthogonal to the face \citep{moukalled2016finite}.
However, when structured or unstructured non-orthogonal meshes are
employed, the described discretization process loses accuracy since
the surface normal vector $\mathbf{s}_{f}$ and $\mathbf{r}$ are
no longer collinear, as illustrated in Figure \ref{fig:non-orthogonal-mesh}.
This implies that, for non-orthogonal meshes, $\hat{\mathbf{n}}\cdot\text{\ensuremath{\nabla}}\mathbf{A}_{f}\neq\hat{\mathbf{r}}\cdot\text{\ensuremath{\nabla}}\mathbf{A}_{f}$.
As Equation \eqref{eq:rGrad} yields $\hat{\mathbf{r}}\cdot\text{\ensuremath{\nabla}}\mathbf{A}_{f}$,
and to evaluate the Laplacian we require computing $\hat{\mathbf{n}}\cdot\text{\ensuremath{\nabla}}\mathbf{A}_{f}$,
an alternative scheme must be employed.

Although it is not possible to calculate the face normal gradient
implicitly, we can calculate the full gradient of the magnetic vector
potential at the cell center explicitly as
\begin{equation}
\overline{\nabla\mathbf{A}_{C}}=\frac{1}{V}\sum\mathbf{s}_{f}\cdot\mathbf{A}_{f}.\label{celGrad}
\end{equation}

This expression allows finding the face gradient explicitly using
linear interpolation in the form
\begin{equation}
\overline{\nabla\mathbf{A}_{f}}=w\text{\ensuremath{\nabla}}\mathbf{A}_{C}+\left(1-w\right)\text{\ensuremath{\nabla}}\mathbf{A}_{E},\label{InterpFaceGrad}
\end{equation}
being $w$ an interpolation factor, and the overbar denotes values
that are calculated explicitly.

Now, we split the evaluation of the face gradient in a orthogonal
component $\text{\ensuremath{\nabla}}\mathbf{A}_{f}^{o}$ and a non-orthogonal
component $\text{\ensuremath{\nabla}}\mathbf{A}_{f}^{n}$ as
\begin{equation}
\text{\ensuremath{\nabla}}\mathbf{A}_{f}=\text{\ensuremath{\nabla}}\mathbf{A}_{f}^{o}+\text{\ensuremath{\nabla}}\mathbf{A}_{f}^{n}.\label{faceGrad}
\end{equation}

We can then extract the normal direction of each gradient to obtain
\begin{equation}
\begin{aligned}\text{\ensuremath{\hat{\mathbf{n}}}\ensuremath{\mathbf{\cdot}\nabla}}\mathbf{A}_{f} & =\hat{\mathbf{n}}\mathbf{\cdot}\text{\ensuremath{\nabla}}\mathbf{A}_{f}^{o}+\hat{\mathbf{n}}\mathbf{\cdot}\text{\ensuremath{\nabla}}\mathbf{A}_{f}^{n}\end{aligned}
.\label{gradSplit}
\end{equation}

Next, we can approximate both the orthogonal component and the LHS
of Eq. \eqref{gradSplit} using the explicit gradient in Eq. \eqref{InterpFaceGrad}
as $\hat{\mathbf{n}}\mathbf{\cdot}\text{\ensuremath{\nabla}}\mathbf{A}_{f}^{o}\simeq\hat{\mathbf{r}}\mathbf{\cdot}\overline{\nabla\mathbf{A}_{f}}$
and $\text{\ensuremath{\hat{\mathbf{n}}}\ensuremath{\mathbf{\cdot}\nabla}}\mathbf{A}_{f}\simeq\text{\ensuremath{\hat{\mathbf{n}}}\ensuremath{\mathbf{\cdot}}}\overline{\nabla\mathbf{A}_{f}}$.
Therefore, we can approximate the non-orthogonal component of the
face normal gradient in Eq. \eqref{InterpFaceGrad} explicitly as
\begin{equation}
\hat{\mathbf{n}}\mathbf{\cdot}\overline{\nabla\mathbf{A}_{f}^{n}}=\left(\hat{\mathbf{n}}-\hat{\mathbf{r}}\right)\mathbf{\cdot}\overline{\nabla\mathbf{A}_{f}}.\label{nOrthGrad}
\end{equation}

Finally, injecting Eqs. \eqref{nOrthGrad} and \eqref{eq:rGrad} into
Eq. \eqref{gradSplit} allows us to find an implicit-explicit expression
for the face normal gradient of the magnetic vector potential, namely
\[
\text{\ensuremath{\hat{\mathbf{n}}}\ensuremath{\mathbf{\cdot}\nabla}}\mathbf{A}_{f}=\frac{\mathbf{A}_{C}-\mathbf{A}_{E}}{\left|\mathbf{r}\right|}+\left(\hat{\mathbf{n}}-\hat{\mathbf{r}}\right)\mathbf{\cdot}\overline{\nabla\mathbf{A}_{f}}.
\]

This expression must be iterated until convergence in a loop that,
in the context of fluid mechanics, is often called non-orthogonal
correction.

\subsubsection{Induced and bound current discretization}

The second term in the RHS of Equation \eqref{eq:dis1} is calculated
explicitly using a Gauss scheme; this yields
\begin{equation}
\intop_{V_{i}}\text{\text{\ensuremath{\nabla\cdot}}}\left(\chi\text{\ensuremath{\nabla}}\widetilde{\mathbf{A}}\right)\,\text{d}V_{i}\approx\sum\chi_{f}\mathbf{s}_{f}\mathbf{\cdot}\text{\ensuremath{\nabla}}\widetilde{\mathbf{A}}_{f},\label{eq:disc-MIi}
\end{equation}

where $\chi_{f}$ is the magnetic susceptibility at face $f$. The
susceptibility at internal faces is calculated using linear interpolation
as
\begin{equation}
\chi_{f}=\left(1-w\right)\chi_{C}+w\chi_{E},
\end{equation}

being $w=\nicefrac{\left|\mathbf{r}_{Cf}\right|}{\left|\mathbf{r}_{EC}\right|}$
the interpolation factor.

At interface faces, the material susceptibility is likely to have
discrete jumps; using interpolation at this faces would smear the
susceptibility field unrealistically. Therefore, we evaluate the face
susceptibility as 
\begin{equation}
\chi_{f}=\chi_{C}+\mathbf{r}_{Cf}\mathbf{\cdot}\nabla\chi_{C},
\end{equation}
for which the susceptibility cell gradient is evaluated explicitly
using a Gauss scheme, 
\begin{equation}
\nabla\chi_{C}=\frac{1}{V_{C}}\sum\chi_{f}\mathbf{s}_{f}.
\end{equation}

To complete the derivation of the discretization scheme, we discretize
the bound current explicitly using the Hodge dual $\text{*}\left(\right)$;
this yields

\begin{equation}
\intop_{V_{i}}\text{\ensuremath{\nabla\times}}\mathbf{M}\,\text{d}V\approx2\text{*}\left(\nabla\widetilde{\mathbf{M}}_{C}\right).
\end{equation}

This method allows choosing a convenient scheme for the gradient of
the magnetization.

\subsection{Solution scheme\label{subsec:Under-relaxation-of-the}}

The numerical discretization presented in the last section results
in coupled sets of nonlinear algebraic equations involving variables
corresponding to different media. In order to solve this multi-domain
problem, we can opt for a monolithic or a partitioned approach. The
monolithic or direct approach addresses all domains simultaneously,
solving a single set of algebraic equations encompassing all relevant
variables. On the other hand, the partitioned or segregated approach
solves the system of equation of each domain individually and sequentially,
treating the coupling terms associated to the interactions between
regions explicitly \citep{mehl2016parallel,heil2008solvers}.

The monolithic approach generates variables that exhibit heterogeneity
as they represent discretization across diverse domains; thus, scalability,
parallelism, preconditioning and performance are commonly an issue
\citep{heil2008solvers}. The monolithic solution of Eq. \eqref{eq:dis1}
with the embedded discontinuity interface conditions in Eq. \eqref{eq:bc-7}
poses a question about stability issues in strongly coupled problems.

In a segregated multi\textendash region scheme solving of a large
system is avoided at the expense of iterations between smaller systems;
typically, the main concerns are stability and implementation \citep{degroote2010stability}.
Also, field transfer between domains, convergence acceleration and
synchronization of region solutions are numerical and computational
difficulties that must be addressed. However, we argue that the partitioned
approach is still more appealing because existing field operation
libraries, such as OpenFOAM \citep{weller1998tensorial}, are equipped
with several discretization, mapping and convergence checking algorithms
that can be straightforwardly combined to develop a effective segregated
solver. In this paper, this is the approach we adopt.

The solution scheme we propose is as follows: suppose a two\textendash region
magnetic system, then the system matrix corresponding to the discretized
version of Eq. \eqref{eq:A-15} is

\begin{equation}
\left[\begin{array}{cc}
\mathbf{L}_{11}+\mathbf{N}_{11}\left(\mathbf{A}_{1}\right) & \mathbf{N}_{12}\left(\mathbf{A}_{2}\right)\\
\mathbf{N}_{21}\left(\mathbf{A}_{1}\right) & \mathbf{L}_{22}+\mathbf{N}_{22}\left(\mathbf{A}_{2}\right)
\end{array}\right]\left[\begin{array}{c}
\mathbf{A}_{1}\\
\mathbf{A}_{2}
\end{array}\right]=\left[\begin{array}{c}
\mathbf{F}_{1}\\
\mathbf{F}_{2}
\end{array}\right],\label{eq:SM}
\end{equation}
being $\mathbf{L}_{jj}$ are linear operators, $\mathbf{N}_{jj}$
are nonlinear operators, $\mathbf{A}_{1}$ and $\mathbf{A}_{2}$ are
the vector potentials for each region, and $\mathbf{F}_{1}$ and $\mathbf{F}_{2}$
correspond to the explicit discretization of the RHS of Eq. \eqref{eq:A-15}.

A direct solution of Eq. \eqref{eq:SM} can be obtained using linearization
schemes available for the solution of nonlinear problems, such as
the Newton\textendash Raphson method, the Picard method, etc. Instead,
we choose to use a block iterative scheme of the form
\begin{equation}
\mathbf{L}_{11}\mathbf{A}_{1}^{i}=\mathbf{F}_{1}-\mathbf{N}_{11}\left(\mathbf{A}_{1}^{i-1}\right)\mathbf{A}_{1}^{i-1}-\mathbf{N}_{12}\left(\mathbf{A}_{2}^{i-1}\right)\mathbf{A}_{2}^{i-1},
\end{equation}

and

\begin{equation}
\mathbf{L}_{22}\mathbf{A}_{2}^{i}=\mathbf{F}_{2}-\mathbf{N}_{22}\left(\mathbf{A}_{2}^{i-1}\right)\mathbf{A}_{2}^{i-1}-\mathbf{N}_{21}\left(\mathbf{A}_{1}^{i}\right)\mathbf{A}_{1}^{i},\label{eq:Algorithm2}
\end{equation}

where $i$ is the iteration counter index. The use of the last available
solution for region 1 in the last term in the RHS of Eq. \eqref{eq:Algorithm2}
indicates that this iteration loop is a block Gauss\textendash Seidel
scheme (BGS).

This algorithm is effective for solving the coupled multi\textendash region
systems; however, the convergence is greatly affected by the sequential
nature of the solver. With highly non\textendash orthogonal or skewed
grids and ill conditioned systems, large variations in $\mathbf{A}$
occur during the solution process, which can cause loss of convergence.
To promote convergence and stabilize the iterative process, a practical
approach is to under\textendash relax certain fields between iterations.

Because the RHS of a region's system is a function of second derivatives
of $\mathbf{A}$, relaxing $\mathbf{A}$ does not always stabilize
the solution. It is not uncommon to find that for cases involving
highly permeable materials the solution diverges. Also, as Eq. \ref{eq:bc-7}
shows, the interface boundary condition is a function of a jump in
the magnetic field at both sides of the interface. Since the BGS scheme
solves each region at a time, during the region sequential loop, one
of the interface neighboring cell center values is updated, while
the others do not. This often generates large oscillations in the
interface jumps.

In the numerical implementation, the interface boundary condition
feeds the Laplacian operator with the face normal gradient of the
vector potential, 
\begin{equation}
\frac{\partial\mathbf{A}_{e}^{C}}{\partial e_{n}}=a_{C}\mathbf{A}^{C}+a_{E}\mathbf{A}^{E}+a_{K}\mathbf{K},\label{dAn}
\end{equation}
being $a_{C}$ the gradient internal factor and $a_{E}$ and $a_{K}$
gradient boundary factors \citep{moukalled2016finite} and the generalized
surface current given by \ref{eq:bc4}.

Since the free surface current $\mathbf{K}_{f}$ and the permanent
magnetization $\mathbf{M}$ are generally constant and uniformly distributed,
the oscillation of the surface current is generated by $\Delta\left(\chi\mathbf{B}\right)$;
this suggests that its relaxation would certainly help convergence.
The first term in the RHS of Eq. \eqref{dAn} adds to the diagonal
of the system matrix, while the other two add to the RHS. For the
case of interface faces, the second term is not a matrix term because
the neighboring cell center value belongs in a different region. In
regards to the third term, it must be evaluated explicitly because
the interface current cannot be written in terms of $\mathbf{A}^{C}$
since the curl couples the vector components. Therefore, both the
second and third terms may generate sharp changes in the RHS of the
region system, leading to instabilities in poor conditioned systems.

In the same direction, the divergence term in the RHS of Eq. \eqref{eq:A-15}
is discretized explicitly as given by Equation \ref{eq:disc-MIi};
implicit integration is not possible since the skew operator couples
the vector components of the vector potential. In this circumstances,
the vector potential gradient can generate values with an nonphysically
large magnitude in highly non-orthogonal meshes. A large gradient
value in one cell can contribute sufficiently to the source term of
an equation to cause unboundedness in the solution.

In order to alleviate the mentioned issues, we use field relaxation.
The main objective of implementing relaxation is to slow down the
changes of the solution variable during the iteration process, which
helps to avoid divergence due to mesh non-orthogonality.

Experimental findings suggest that both implicit and explicit relaxation
can improve the convergence of the proposed solver. For the case of
explicit relaxation, a field $\phi$ in the current iteration is evaluated
in terms of the actual solution of the region system $\phi_{\mathrm{c}}$,
and the previous iteration solution $\phi_{\mathrm{o}}$ as \citep{moukalled2016finite}
\begin{equation}
\phi=\phi_{\mathrm{o}}+\lambda\left(\phi_{\mathrm{c}}-\phi_{\mathrm{o}}\right),\label{eq:explicit-relaxation-1}
\end{equation}
where $\lambda$ represents the relaxation factor. 

As it can be seen from the expression above, the current field $\phi$
is not the calculated $\phi_{\mathrm{c}}$ but a fraction of it, this
reduces oscillations for values of $\lambda<1$. When convergence
is reached, the correction tends to zero and $\phi_{\mathrm{c}}=\phi_{\mathrm{o}}$.
Although reducing the value of $\lambda$ increases stability, the
correction is also reduced and the convergence is slowed down.

Our experimental findings show that when explicit relaxation is used,
the proposed iterative block algorithm diverges for most high permeability
configurations. Under this circumstances, the following implicit relaxation
strategy has shown to improve stability. 

Consider the matrix equation 
\begin{equation}
\mathbf{C}\,\boldsymbol{\phi}=\mathbf{D},\label{eq:generuc-matrix-equation}
\end{equation}
being $\mathbf{C}$ a coefficients matrix of size $n\times n$, $\boldsymbol{\phi}$
the solution vector and $\mathbf{D}$ a source term. For a single
cell $P$, the equation is
\begin{equation}
C_{P1}\phi_{1}+C_{P2}\phi_{2}+\ldots+C_{P(P-1)}\phi_{P-1}+C_{PP}\phi_{P}+C_{P(P+1)}\phi_{P+1}+\ldots+C_{Pn}\phi_{n}=D_{P}.
\end{equation}

Since we want the solution for cell $P$, i.e. $\phi_{P}$, we can
collect the other terms together, which, using summation notation,
gives
\begin{equation}
C_{P}\phi_{P}+\sum_{nb}C_{nb}\phi_{nb}=D_{P},\label{eq:implicit-relaxation-1}
\end{equation}

where $nb$ stands for neighbor cells. Thus, we can express the solution
for cell $P$ in the form
\begin{equation}
\phi_{P}=\frac{1}{C_{P}}\left(-\sum_{nb}C_{nb}\phi_{nb}+D_{P}\right).
\end{equation}

This value corresponds to the component $P$ of $\boldsymbol{\phi}$.
Now, explicitly relaxing $\phi_{P}$
\begin{equation}
\phi_{P}=\phi_{\mathrm{o}}+\lambda\left[\frac{1}{C_{P}}\left(-\sum_{nb}C_{nb}\phi_{nb}+D_{P}\right)-\phi_{\mathrm{o}}\right].
\end{equation}
Multiplying both sides by $\nicefrac{C_{P}}{\lambda}$ and rearranging
terms results in
\begin{equation}
\frac{C_{P}}{\lambda}\phi_{P}+\sum_{nb}C_{nb}\phi_{nb}=D_{P}+\left(\frac{1-\lambda}{\lambda}\right)C_{P}\phi_{\mathrm{o}}.\label{eq:implicit-relaxation-2}
\end{equation}

Writing Equation \ref{eq:implicit-relaxation-2} in matrix form yields
\begin{equation}
\left[\begin{array}{ccc}
\frac{C_{11}}{\lambda} & \text{\ensuremath{\ldots}} & C_{1n}\\
\vdots & \ddots & \vdots\\
C_{n1} & \ldots & \frac{C_{nn}}{\lambda}
\end{array}\right]\left[\begin{array}{c}
\phi_{1}\\
\vdots\\
\phi_{n}
\end{array}\right]=\left[\begin{array}{c}
D_{1}+\left(\frac{1-\lambda}{\lambda}\right)C_{11}\phi_{1,\mathrm{o}}\\
\vdots\\
D_{n}+\left(\frac{1-\lambda}{\lambda}\right)C_{nn}\phi_{n,\mathrm{o}}
\end{array}\right].
\end{equation}

From the equation above, it is clear that implicit relaxation enhances
the diagonal dominance of the system matrix, i.e. each cell's contribution
to its own solution becomes more significant compared to the influence
of neighboring cells. The increased diagonal dominance promotes stability,
as it ensures smoother propagation of changes throughout the system,
thereby reducing the likelihood of oscillations or divergence during
the iterative process. Another consequence of implicit relaxation
is an increase in the source term, which typically slows down convergence. 

At this point, it is important to remark that there are two primary
strategies for implementing under-relaxation at the multi-region solver
level. The first involves a region-wise application, while the second
applies simultaneous relaxation to all regions. The former approach
introduces a significant dependency on the order in which regions
are solved, as an under-relaxed solution may propagate to neighboring
regions through the interface boundary conditions. Conversely, such
concerns are mitigated in the latter approach, as the relaxation of
the right-hand side of the equilibrium equation is uniformly applied
across all regions before commencing the iterative solution loop.

\section{Numerical experiments\label{sec:Numerical-experiments}}

In this section, we present an assessment of the method's performance.
We evaluate its capability to accurately capture discontinuous distributions
of magnetic fields in configurations comprising permeable, permanently
magnetized, and current-carrying media by comparing its results with
those obtained using the Finite Element Method. All test cases are
two-dimensional and are constructed using triangular non-structured
meshes.

\subsection{Magnetized and permeable media interaction - Case 1}

\begin{figure}[H]
\begin{centering}
\includegraphics{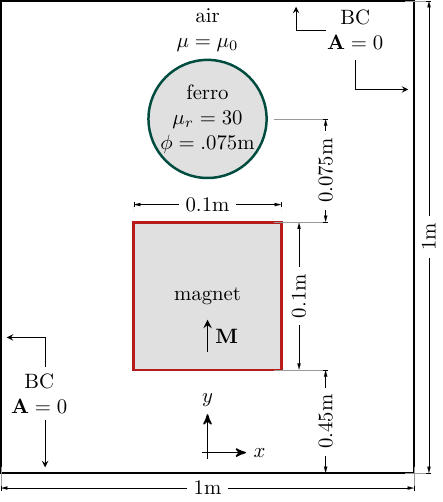}
\par\end{centering}
\caption{Schematic of benchmark Case 1. \label{fig:test-1}}
\end{figure}

The first test consists of a squared magnet interacting with a permeable
material. As Figure \ref{fig:test-1} shows, the configuration consists
of a cylindrical shaped permeable material, called \emph{ferro,} with
diameter $\phi=0.075\,\text{m}$ and a squared permanently magnetized
material, called \emph{magnet}, with sides of length $0.1\,\text{m}$.
The ferro's relative permeability is $\mu_{r}=30$ and the permanent
magnetization of the magnet is $\mathbf{M}=9.75\times 10^{5} \ \text{A m$^{-1}$}\mathbf{e}_{y}$;
being $\mathbf{e}_{y}$ the unit vector pointing to vertical direction.
Both the magnet and the ferro are immersed in an air volume of dimensions
$1\,\text{m}\times1\,\text{m}\times1\,\text{m}$. Air magnetic properties
are assumed to be $\mu_{r}=1$ and $\mathbf{M}=\mathbf{0}$. At the
outer boundaries of the air volume we specify the boundary condition
$\mathbf{A}=\mathbf{0}$; while at all interfaces we specify the boundary
condition in Eq. \eqref{dAn}.

\begin{figure}[H]
\begin{centering}
\includegraphics[width=80mm]{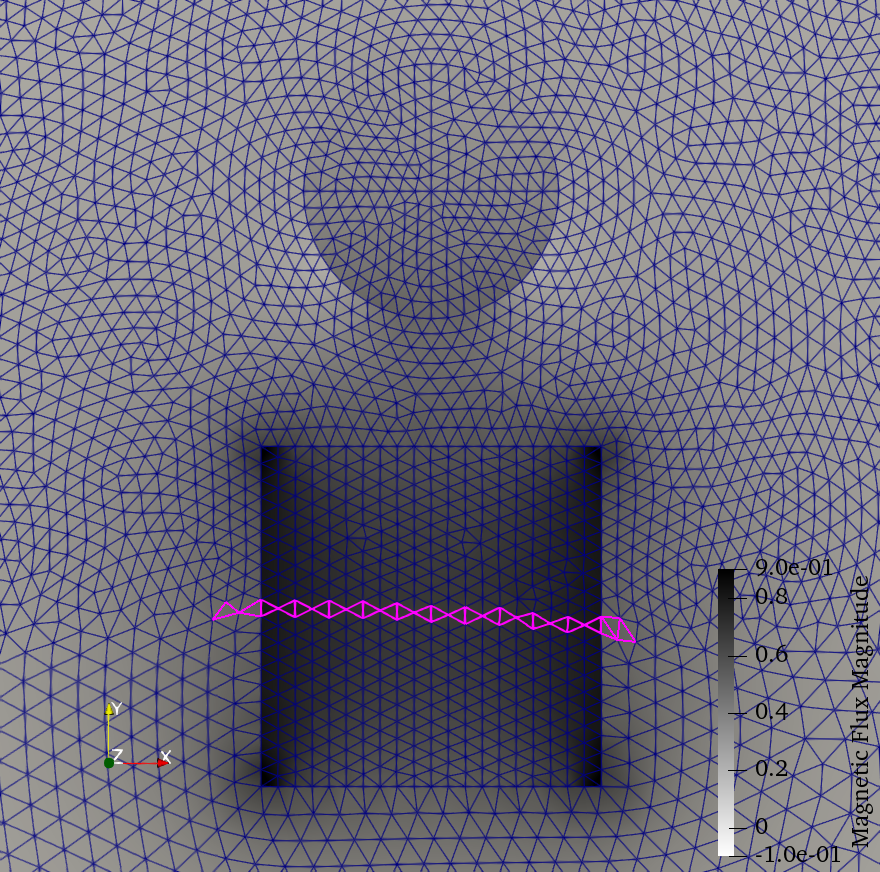}
\par\end{centering}
\caption{Benchmark Case 1 mesh. \label{fig:mesh-1}}
\end{figure}

The mesh density around the ferro and magnet has fixed edge size of
$5\ times10^{-3}\,\text{m}$; Figure \ref{fig:mesh-1} shows the mesh.

\begin{figure}[H]
\begin{centering}
\includegraphics[width=80mm]{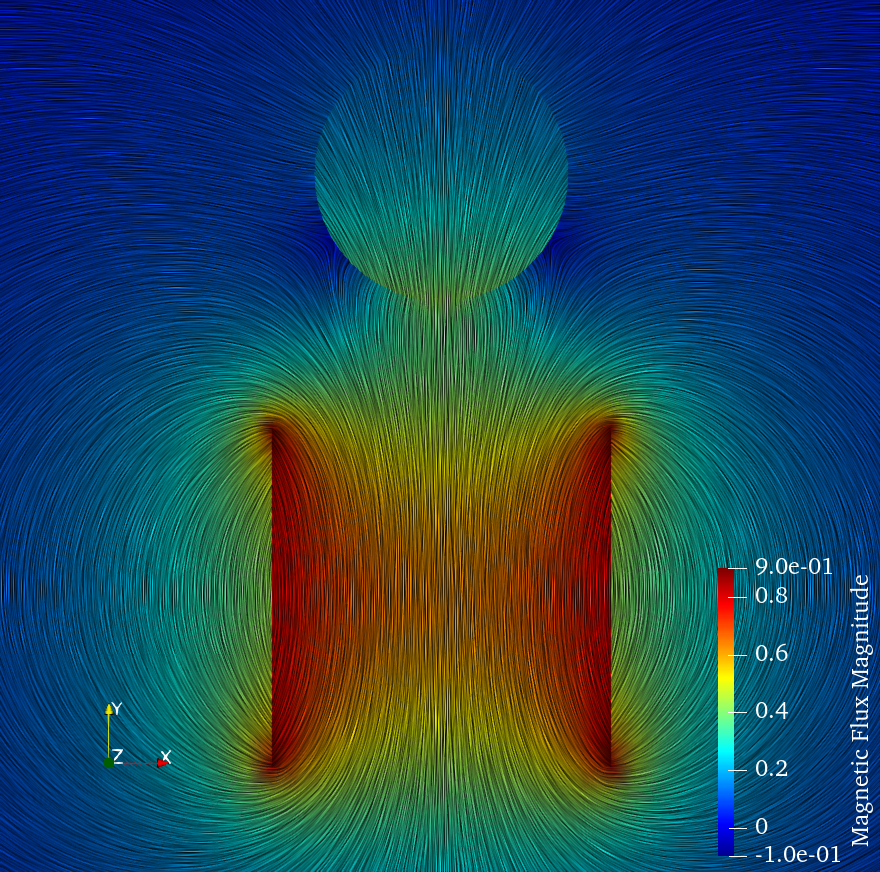}
\par\end{centering}
\caption{Magnetic field SLIC for Case 1 \textendash{} FEM solution. \label{fig:lic-1}}
\end{figure}

\begin{figure}[H]
\begin{centering}
\includegraphics[width=80mm]{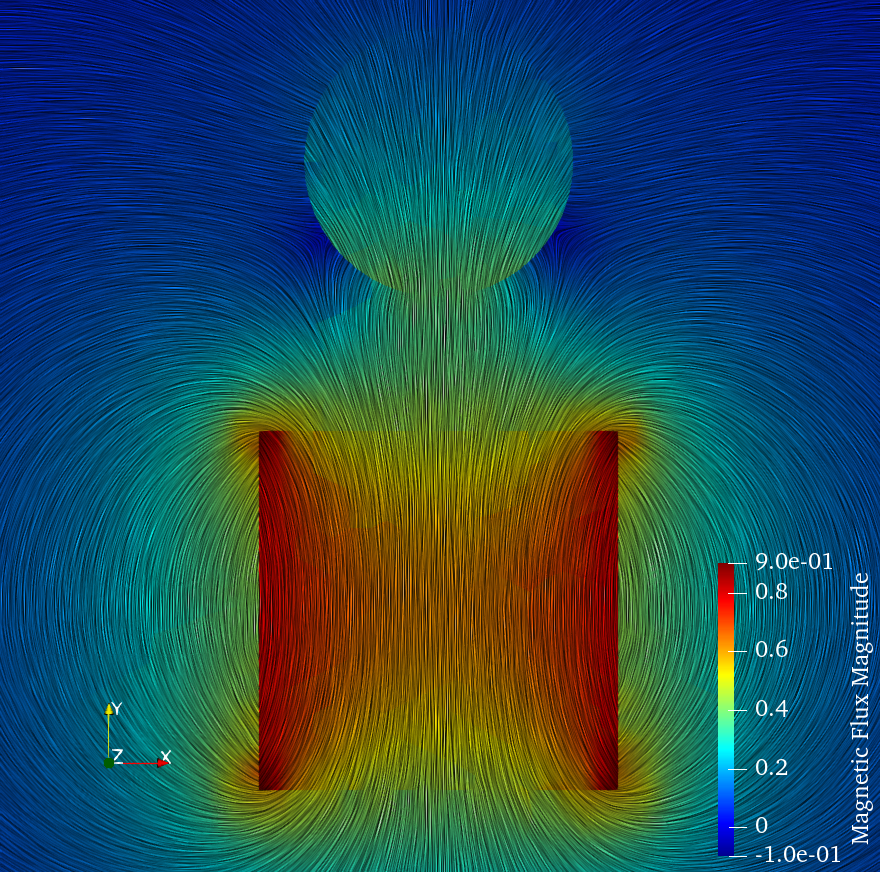}
\par\end{centering}
\caption{Magnetic field SLIC for Case 1 \textendash{} Present FVM solution.
\label{fig:lic-2}}
\end{figure}

Figures \ref{fig:lic-1} and \ref{fig:lic-2} show the surface line
integral convolution (SLIC) of the magnetic field obtained with the
present formulation and the FEM. From a qualitative viewpoint, it
is clear that both approaches agree very well. The discontinuities
of the tangential component of the magnetic at the interfaces between
both the magnet and the permeable material and air between different
media are very similarly captured by both approaches.

\begin{figure}[H]
\begin{centering}
\includegraphics{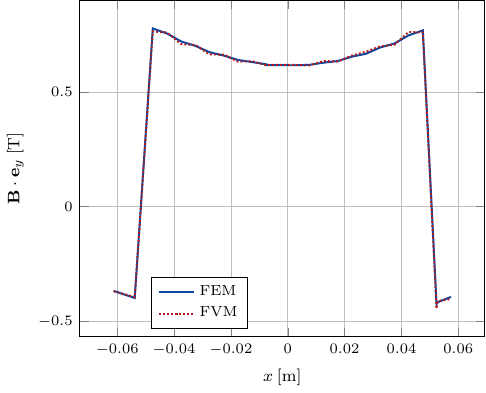}
\par\end{centering}
\caption{Vertical component of the magnetic field.\label{fig:quali-1}}
\end{figure}

\textcolor{black}{Figure \ref{fig:quali-1} shows the change in the
$y$-component of the magnetic field at the cell centroids through
the arbitrary path shown in Figure} \ref{fig:mesh-1}. It is important
to note that since the Finite Element formulation gives magnetic field
values at the element edges, in order to compare the results against
the present cell-centered FVM, it was necessary to interpolate the
FEM data to the element centroid. As it can be seen, the formulation
is able to capture correctly the discontinuous distribution of the
magnetic field.

\subsection{Magnetized and permeable media interaction - Case 2}

\begin{figure}[H]
\begin{centering}
\includegraphics{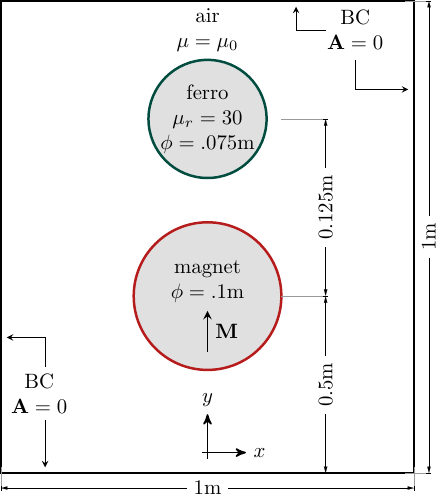}
\par\end{centering}
\caption{Schematic of benchmark Case 2. \label{fig:test-2}}
\end{figure}

Expanding upon the initial comparison, we tested the FV method with
the configuration shown in Figure \ref{fig:test-2}, which features
a domain composed of both cylindrical shaped magnet and permeable
materials, with diameters $\phi=0.1\,\text{m}$ and $\phi=0.075\,\text{m}$
respectively. The magnetization is kept as $\mathbf{M}=9.75\times 10^{5}\,\text{A $m^{-1}$}\mathbf{e}_{y}$.

\begin{figure}[H]
\begin{centering}
\includegraphics[width=80mm]{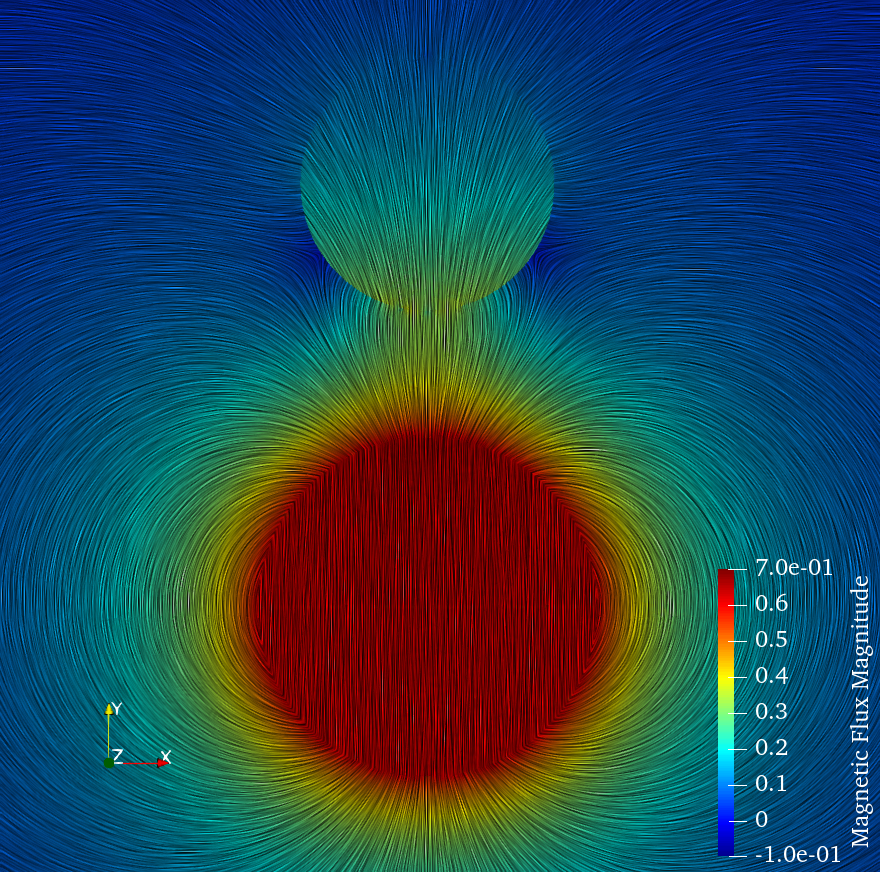}
\par\end{centering}
\caption{Magnetic field SLIC. \label{fig:lic-3}}
\end{figure}

The results of the numerical experiment are shown in Figure \ref{fig:lic-3}.
Qualitatively, the discontinuities in direction and magnitude of $\mathbf{B}$
arising at the interfaces are captured correctly by the solver. Then,
far from the ferro-air and magnet-air interfaces, the results are
consistent with \citep{saravia2021finite} and \citep{riedinger2023single}
and the FEM edge prediction shown in the previous section.

\begin{figure}[H]
\begin{centering}
\includegraphics{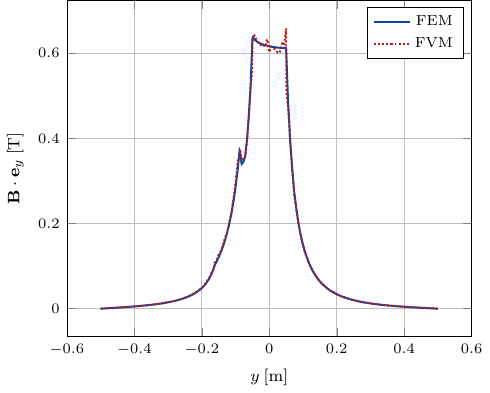}
\par\end{centering}
\caption{Vertical component of the magnetic field.\label{fig:quanti-2}}
\end{figure}

This agreement is confirmed by the quantitative analysis presented
in Figure \ref{fig:quanti-2}, where the variation in the $y$-component
of the magnetic field is plotted along a vertical path through the
center of the domain.

\subsection{Magnetized, permeable and volumetric current media interaction -
Case 3\label{subsec:Volumetric-current-with}}

\begin{figure}[H]
\begin{centering}
\includegraphics{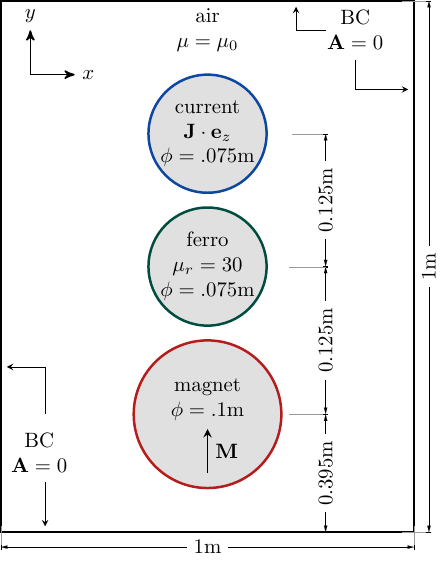}
\par\end{centering}
\caption{Benchmark test schematic. \label{fig:test-3}}
\end{figure}

Next, we assess the method's behavior in a case with strong interaction
among three curved materials. The configuration is shown in Figure
\ref{fig:test-3}, with the benchmark test being composed of a magnet
and ferro materials with the same physical characteristics as in Cases
1 and 2, and additionally a current carrying media of diameter $\phi=0.075\,\text{m}$
and current vector $\mathbf{J}=2.5\times 10^{7} \ \text{A $m^{-2}$} \ \mathbf{e}_{z}$.

\begin{figure}[H]
\begin{centering}
\includegraphics[width=80mm]{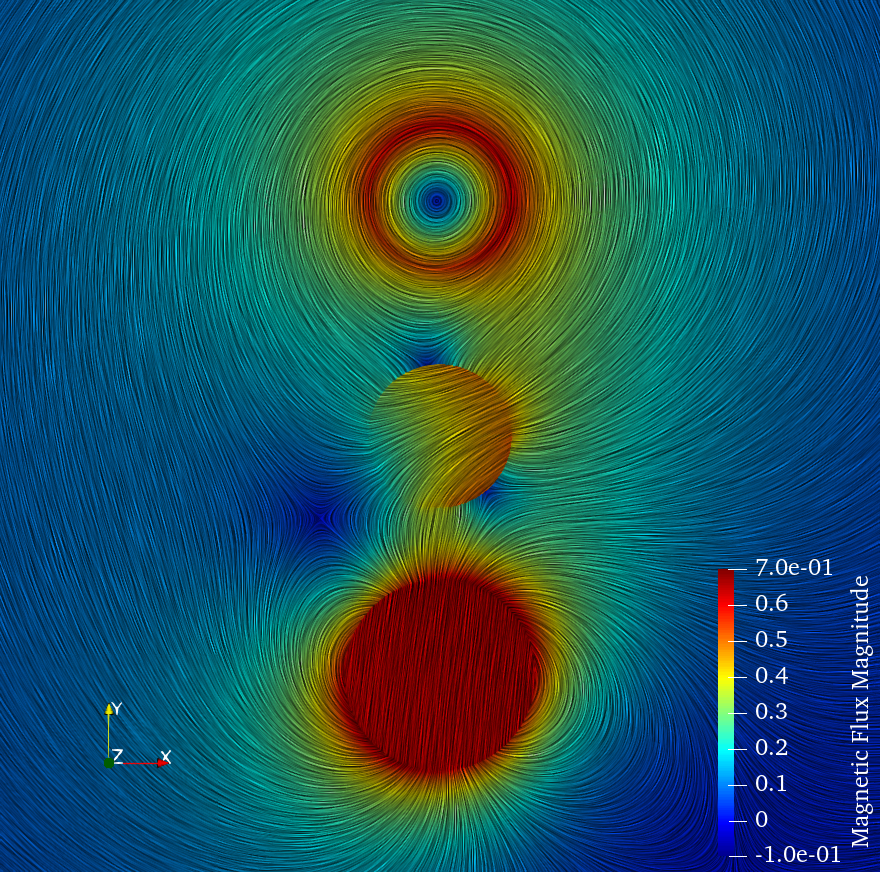}
\par\end{centering}
\caption{Magnetic field SLIC. \label{fig:lic-4}}
\end{figure}

The resulting magnetic field may be observed in Figure \ref{fig:lic-4}.
The current carrying material is characterized by the lack of surface
currents $\mathbf{K}_{f}$, and thus the RHS of the boundary condition,
Eq. \ref{eq:bc-7}, is null since $\Delta\chi\mathbf{B}=\Delta\mathbf{M}$.
Accordingly, the normal and tangential components of the vector potential
in the interface of the current material are continuous. As it can
be seen qualitatively in the figure, the method is able to capture
the changes in the direction and magnitude of the vector magnetic
field in all three materials.

\begin{figure}[H]
\begin{centering}
\subfloat[Vertical component of the magnetic field.]{\begin{centering}
\includegraphics{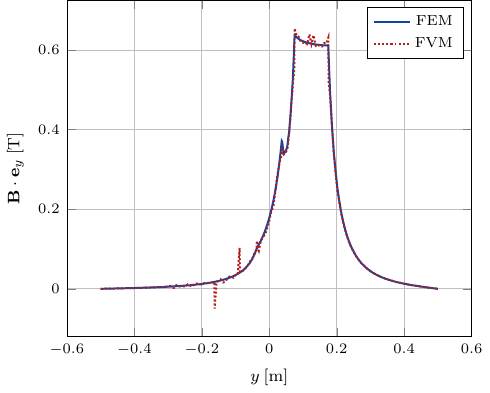}
\par\end{centering}
}\subfloat[Horizontal component of the magnetic field.]{\begin{centering}
\includegraphics{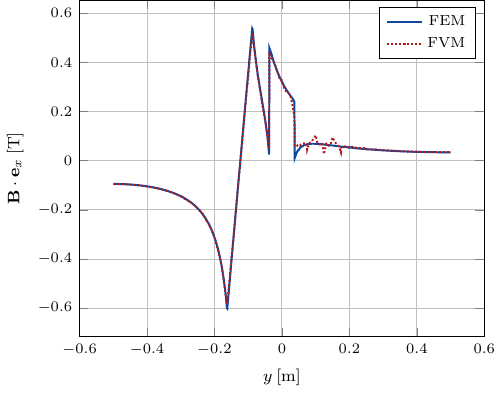}
\par\end{centering}
}
\par\end{centering}
\begin{centering}
\subfloat[Depth component of the magnetic vector potential.]{\begin{centering}
\includegraphics{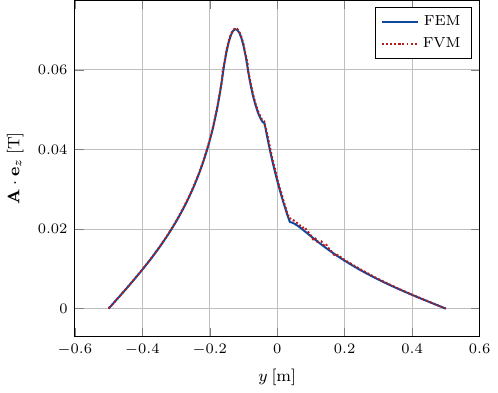}
\par\end{centering}
}
\par\end{centering}
\caption{Quantitative assessment of Case 3.\label{fig:quanti-3}}
\end{figure}

Quantitatively, the results are illustrated in Figure \ref{fig:quanti-3}.
In this case, the changes in the $x$- and $y$-components of the
magnetic field, as well as the change in the $z$-component of the
magnetic vector potential, are plotted. The method successfully captures
the variations in the magnetic field despite the strong interactions
between the media. Notably, the plot of the magnetic vector potential
closely aligns with the FEM results; whereas the magnetic field plots
exhibit more error due to the additional step of taking the curl of
the solution $\mathbf{A}$.

\subsection{Non-orthogonality test - Case 4}

\begin{figure}[H]
\begin{centering}
\includegraphics{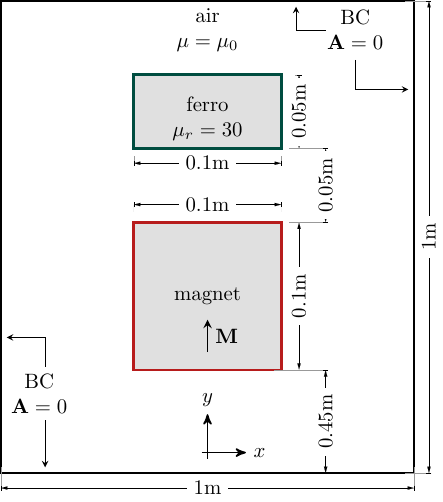}
\par\end{centering}
\caption{Case 4 schematic. \label{fig:test-4}}
\end{figure}

Finally, we solve the configuration depicted in Figure \ref{fig:test-4},
which is taken from \citep{saravia2021finite,riedinger2023single}.
The objective of this test is to compare the results obtained with
orthogonal and non-orthogonal meshes. The configuration consists of
two squared materials: a magnet of sides $0.1\,\text{m}$ and a ferro
of size $0.1\,\text{m}\times0.05\,\text{m}$, surrounded by a volume
of air of dimensions $1\,\text{m}\times1\,\text{m}$.

\begin{figure}[H]
\begin{centering}
\includegraphics[width=80mm]{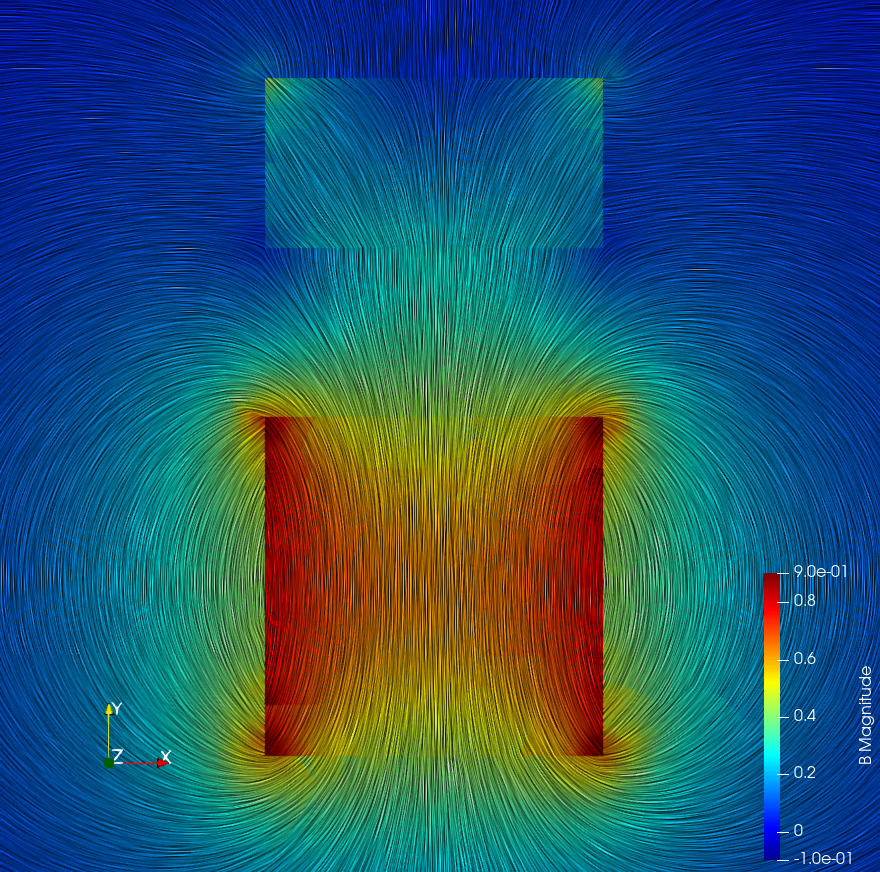}
\par\end{centering}
\caption{Magnetic field SLIC \textendash{} Non-orthogonal mesh. \label{fig:lic-5}}
\end{figure}

\begin{figure}[H]
\begin{centering}
\includegraphics[width=80mm]{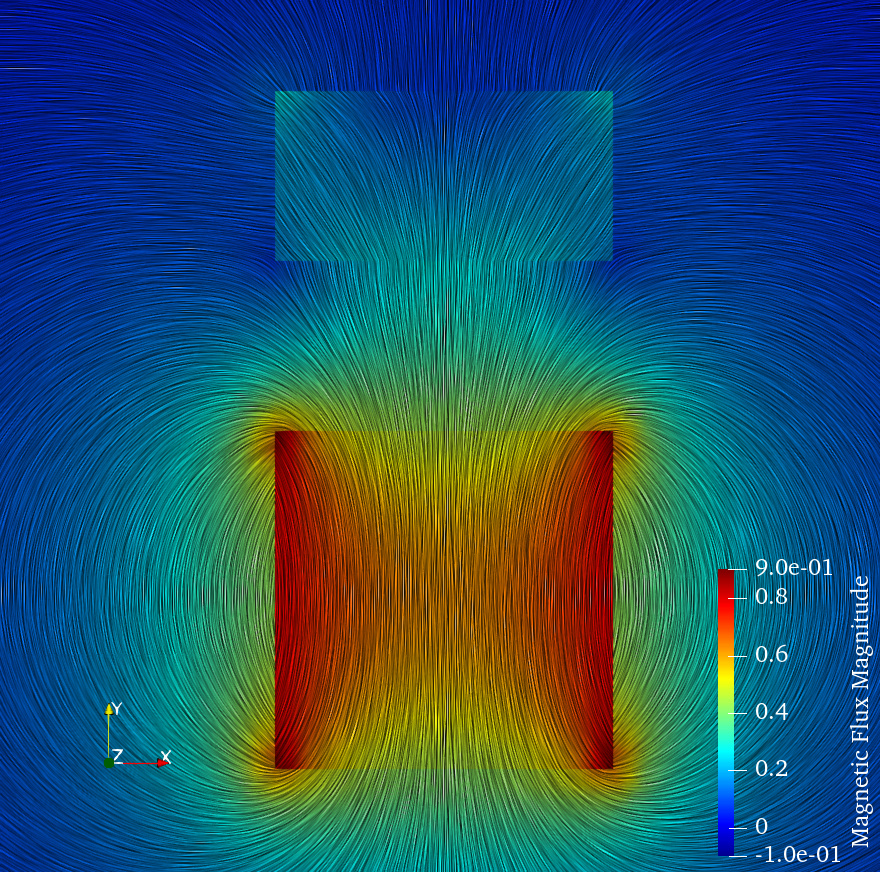}
\par\end{centering}
\caption{Magnetic field SLIC \textendash{} Orthogonal mesh. \label{fig:lic-6}}
\end{figure}

The results of this test are shown in Figure \ref{fig:lic-5}, while
Figure \ref{fig:lic-6} depicts the same results but with an orthogonal
and uniform mesh, maintaining similar mesh density. As it can be seen
from the qualitative assessment, the calculation of the magnetic field
coincides regardless of the mesh. However, in the non-orthogonal case,
convergence is only achieved taking into account considerations such
as field relaxation. This is not case for the orthogonal experiment,
where these considerations are not necessary and therefore not implemented.

The results of this test are presented in Figure \ref{fig:lic-5}
and \ref{fig:lic-6}. The figures show the magnetic field ditribution
obtained with non-orthogonal and orthogonal meshes, respectively.
Both meshes have comparable density. The qualitative assessment indicates
that regardless of the mesh topology, the algorithm is able to predict
the same magnetic field. A key point is that the convergence of the
non-orthogonal case requires the use of a relaxation technique. In
this case, we have used implicit relxation with a relaxation factor
$\lambda=0.8$; larger values of $\lambda$ cause divergence.

As a general observation, we note that behavior of both Gauss-linear
and least-squares gradient methods for the evaluation of the non-orthogonal
face normal gradient in Eq. \eqref{nOrthGrad} were evaluated. It
was found that the least-squares gradient typically yields divergent
solutions, whereas the Gauss gradient method consistently demonstrates
stability.

Numerical experiments also indicate that convergence is governed by
the bound current term $\nabla\cdot\left(\chi\nabla\widetilde{\mathbf{A}}\right)$.
Convergence is often guaranteed by under-relaxing this term with $\lambda\approx0.8$.
The relaxation of the magnetic vector potential and/or the magnetic
field does not seem to help convergence. Comparisons between implicit
and explicit relaxation schemes reveal that implicit relaxation is
generally more efficient, resulting in higher relaxation factors and
faster solutions.

We have also tested the $\chi$ scaling approach suggested in \citep{saravia2021finite}
for the improving the convergence behavior in orthogonal meshes. This
approach does not improve the convergence behavior in non-orthogonal
meshes. 

\section{Conclusions and future work\label{sec:Conclusions-and-future}}

We presented an analysis of discontinuous distributions of magnetic
field within curved media discretized with non-orthogonal grids. A
formulation where the balance equation is based upon the magnetic
vector potential has been developed and tested in diverse numerical
settings with non-orthogonal meshes. Since the interface vector magnetic
potential is provided by the solution rather than being enforced as
an internal boundary condition, numerical experiments have been developed
to asses the magnetic field in the multi-media interface.

One of the key conclusions that have been reached from the numerical
experiments is the considerations neccesary to achieve convergence
with non-orthogonal grids. We have ascertained that field relaxation
is imperative in order to find a solution; precisely, relaxing the
divergence term in the RHS of Equation \ref{eq:A-15} with a factor
$\lambda\sim\mathcal{O}(0.8)$ is needed instead of scaling the parameter
$\chi$ as it has been done in \citep{saravia2021finite} for orthogonal
meshes. This has been verified by comparing the cell centre values
of the Finite Volume framework with interpolated values given by the
FE COMSOL solver.

Once the predictions were verified and contrasted against the FE solver,
further numerical experiments have been benchmarked in order to test
the approach in more complex geometries and with stronger field interactions.
The results suggest that the solver is capable of handling these geometries,
finding solutions regardless of the meshing techniques employed to
discretize the space.

The effects of non-orthogonality have been assessed numerically by
comparing the surface line integral convolution plots with both triangular
and quadratic conformal grids within the same geometry. Although this
comparison is qualitative rather than quantitative, it is useful to
get a general idea of the solver's accuracy with non-orthogonal meshes.

The approach has proven to be effective for the computation of discontinuous
distributions of magnetic field within curved surfaces. Moreover,
the convergence requirements that have been found lay the foundation
for more robust computations of magnetic fields. Thus, it is concluded
that using a FV scheme shall yield the similar results as a FE solver
regardless of the mesh employed.

\appendix

\section{Vector identities}

Being, $c$ a scalar, $\underline{a}$ and $\underline{b}$ vectors,
$\underline{e_{n}}$ a normal unit vector, $\underline{T}$ a tensor,
the following identities hold 
\begin{equation}
\mathrm{\nabla\times}\left(\mathrm{\nabla\times}\,\mathbf{a}\right)=\mathrm{\nabla}\left(\mathrm{\nabla\cdot}\,\mathbf{a}\right)-\mathrm{\nabla^{2}}\,\mathbf{a};\label{eq:I03}
\end{equation}
\begin{equation}
\mathrm{\nabla\times}\left(c\,\mathbf{a}\right)=c\,\mathrm{\nabla\times}\,\mathbf{a}+\mathrm{\nabla}\,c\times\mathbf{a};\label{eq:I04}
\end{equation}
\begin{equation}
\mathbf{\nabla}\mathbf{\cdot}\left(\phi\mathbf{T}\right)\mathbf{=}\phi\mathbf{\nabla}\cdot\mathbf{T}\mathbf{+}\mathbf{T}^{T}\mathbf{\nabla}\phi\label{eq:I09}
\end{equation}
\begin{equation}
\mathrm{\nabla\cdot}\left(\mathrm{\nabla}\,\mathbf{a}\right)\mathbf{=}\mathrm{\nabla^{2}}\,\mathbf{a};\label{eq:I10}
\end{equation}
\begin{equation}
\mathrm{\nabla\cdot}\left(\left(\mathrm{\nabla}\,\mathbf{a}\right)^{T}\right)=\mathrm{\nabla}\left(\mathrm{\nabla\cdot}\,\mathbf{a}\right);\label{eq:I11}
\end{equation}
\begin{equation}
\left(\nabla c\right)\times\left(\nabla\times\mathbf{a}\right)=\left(\nabla\mathbf{a}-\nabla\mathbf{a}^{T}\right)\nabla c.\label{eq:I12}
\end{equation}

\section*{Acknowledgments}

We would like to thank the support from Consejo Nacional de Investigaciones
Científicas y Técnicas and Universidad Tecnológica Nacional through
grant PID-ENUTIBB0007877TC.

\bibliographystyle{elsarticle-num}
\addcontentsline{toc}{section}{\refname}\bibliography{references}

\end{document}